\documentclass[12pt,reqno]{article}

\usepackage[usenames]{color}
\usepackage{amssymb}
\usepackage{graphicx}
\usepackage{amscd}
\usepackage{latexsym}

\usepackage{amsthm}
\newtheorem{theorem}{Theorem}

\newtheorem{proposition}[theorem]{Proposition}

\theoremstyle{definition}

\newtheorem{example}[theorem]{Example}

\newtheorem{conjecture}[theorem]{Conjecture}

\usepackage[colorlinks=true,
linkcolor=webgreen, filecolor=webbrown,
citecolor=webgreen]{hyperref}

\definecolor{webgreen}{rgb}{0,.5,0}
\definecolor{webbrown}{rgb}{.6,0,0}

\usepackage{color}

\usepackage{float}

\usepackage{graphics,amsmath,amssymb}
\usepackage{amsfonts}
\usepackage{latexsym}
\usepackage{epsf}

\begin{document}

\begin{center}
\vskip 1cm{\LARGE\bf Some conjectures on the ratio of Hankel
transforms for sequences and series reversion} \vskip 1cm \large
Paul Barry\\
School of Science\\
Waterford Institute of Technology\\
Ireland\\
\href{mailto:pbarry@wit.ie}{\tt pbarry@wit.ie} \\
\end{center}
\vskip .2 in

\begin{abstract} For each element of certain families of sequences, we study
the term-wise ratios of the Hankel transforms of three sequences
related to that element by series reversion. In each case, the
ratios define well-known sequences, and in one case, we recover the
initial sequence.
\end{abstract}
\section{Introduction} The Hankel transform for sequences (defined below) has attracted an increasing
amount of attention in recent years. The paper \cite{Layman} situated its study within the mainstream of
research into integer sequences, while papers such as \cite{Hankel1} hinted at how the study of
certain Hankel transforms can lead to results concerning classical sequences. That paper exploited a link
between continued fractions and the Hankel transform, as explained by Krattenhaler \cite{Krat}. The best known example of a Hankel
transform for sequences is that of the Catalan numbers. One of the earlier contributors to our stock of
knowledge about the Hankel transform, Christian Radoux, had published several proofs of this result, along with other
interesting examples \cite{Rad1}, \cite{Rad2},\cite{Rad3},\cite{Rad4},\cite{Rad5}. One should also note the interesting umbral interpretation
of the Hankel transform given in \cite{Zeil}. In this paper we indicate that the term-wise ratio of Hankel transforms of
shifted sequences are noteworthy objects of study, giving us more insight into the processes involved in the Hankel transform.

\section{Integer Sequences and Transforms on them}
In this note, we shall consider integer sequences
$$a:\mathbf{N}_0 \to \mathbf{Z} $$ with general term $a_n=a(n)$. Normally, sequences will be
described by their ordinary generating function (o.g.f.), that is, the function $g(x)$ such that
$$g(x)=\sum_{n=0}^{\infty}a_nx^n.$$
We shall study the Hankel transform of sequences in this note. This
is a transformation on the set of integer sequences defined as
follows. Given a sequence $a_n=a(n)$ as described above, we form the
$(n+1)\times(n+1)$ matrix $H_n$ with general term $a(i+j)$, where
$0\le i,j \le n$. Then the Hankel transform $h_n$ of the sequence
$a_n$ is defined by $$ h_n=\det(H_n).$$ Since the elements of the matrix
$H_n$ are elements of an integer sequence, it is clear that $h_n$ is
again an integer sequence. We shall see later that this
transformation is not invertible.

\begin{example} The Catalan numbers $1,1,2,5,14,42\ldots$, defined by $C(n)=\frac{\binom{2n}{n}}{n+1}$ have o.g.f. $\frac{1-\sqrt{1-4x}}{2x}$.
The Hankel transform of the Catalan numbers is the
sequence of all $1$'s. Thus each of the determinants
$$ |1|, \qquad \left|\begin{array}{cc} 1 & 1 \\ 1 & 2 \end{array}\right|,
\qquad \left|\begin{array}{ccc} 1 & 1 & 2\\1 & 2 & 5\\2 & 5 &
14\end{array}\right|, \qquad \ldots$$ has value $1$.  A unique
feature of the Catalan numbers is that the shifted sequence $C(n+1)$
also has a Hankel transform of all $1$'s. An interesting feature of
the Catalan numbers is that the sequence $C(n)-0^n$, or
$0,1,2,5,14,42,\ldots$ has Hankel transform $n$. Of direct relevance
to this note is the fact that the Hankel transform of the sequence
$0,1,1,2,5,14,42,\ldots$ with o.g.f. $\frac{1-\sqrt{1-4x}}{2}$ is
$-n$. This sequence is defined by the series reversion of the
logistic function $x(1-x)$.
\end{example}
\begin{example} The central binomial coefficients $1,2,6,20,70,252,\ldots$, defined by $a_n=\binom{2n}{n}$, have o.g.f. $\frac{1}{\sqrt{1-4x}}$.
The Hankel transform of the central binomial coefficients is given by $h_n=2^n$. That is,
$$ |1|=1, \qquad \left|\begin{array}{cc} 1 & 2 \\ 6 & 20 \end{array}\right|=2,
\qquad \left|\begin{array}{ccc} 1 & 2 & 6\\2 & 6 & 20\\20 & 70 & 252\end{array}\right|=4, \qquad \ldots$$

The sequence $0,1,2,6,20,\ldots$ with o.g.f. $\frac{x}{\sqrt{1-4x}}$
has Hankel transform $-n2^{n-1}$. This is the negative of the
binomial transform (see below) of $n$.
 \end{example}

An important transformation on integer sequences that is invertible
is the so-called \emph{Binomial transform}. Given an integer
sequence $a_n$, this transformation returns the sequence with
general term
$$b_n = \sum_{k=0}^n \binom{n}{k}a_k.$$ If we consider the sequence with general term $a_n$ to be the vector
$\mathbf{a}=(a_0,a_1,\ldots)$
 then we obtain the binomial transform of the sequence by multiplying this (infinite)
 vector by the lower-triangle matrix $\mathbf{B}$ whose $(n,k)$-th element is equal to
$\binom{n}{k}$:
\begin{displaymath}\mathbf{B}=\left(\begin{array}{ccccccc} 1 & 0 & 0
& 0 & 0 & 0 & \ldots \\1 & 1 & 0 & 0 & 0 & 0 & \ldots \\ 1 & 2 & 1
& 0 & 0 & 0 & \ldots \\ 1 & 3 & 3 & 1 & 0 & 0 & \ldots \\ 1 & 4 &
6 & 4 & 1 & 0 & \ldots
\\1 & 5 & 10 & 10 & 5 & 1 &\ldots\\ \vdots & \vdots & \vdots & \vdots & \vdots
& \vdots & \ddots\end{array}\right)\end{displaymath} The inverse transformation
is given by
$$a_n=\sum_{k=0}^n (-1)^{n-k}\binom{n}{k}b_k.$$
If we use ordinary generating functions to describe a sequence, then
the sequence with o.g.f. $g(x)$ will have a binomial transform whose
o.g.f. is given by $\frac{1}{1-x}g(\frac{x}{1-x})$. \newline\newline
It is shown in \cite{Layman} that if $b_n$ is the binomial transform
of the sequence $a_n$, then both sequences have the same Hankel
transform. Thus the Hankel transform is not invertible.
\section{On the series reversion of certain families of generating
functions of sequences} In this note we shall be concerned mainly
with the Hankel transform of sequences whose o.g.f. will be defined
as the series reversion of the o.g.f.'s of certain basic sequences.
Thus in this section, we will briefly recall facts about the
sequences with o.g.f.'s of the forms given by
$\frac{x}{1+{\alpha}x+{\beta}x^2}$, $\frac{x(1-\alpha x)}{1-\beta
x}$ and $x(1-\alpha x)$ as well as their reversions. The first two
families have been studied in \cite{PasTri}.
\begin{example} The family $\frac{x}{1+{\alpha}x+{\beta}x^2}$.
\newline\newline
The sequence with o.g.f. $\frac{x}{1+{\alpha}x+{\beta}x^2}$ has
general term given by $$\sum_{k=0}^{\frac{\lfloor n-1
\rfloor}{2}}\binom{n-k-1}{k}(-\alpha)^{n-2k}(-\beta)^k.$$ The
reversion of the series $\frac{x}{1+{\alpha}x+{\beta}x^2}$, that is,
the solution $u=u(x)$ to the equation
$$\frac{u}{1+\alpha u+\beta u^2}=x$$ is given by $$u(x)=\frac{1-\alpha{x}-\sqrt{1-2\alpha{x}+(\alpha^2-4\beta)x^2}}{2\beta{x}}.$$
The sequence $u_n$ with this o.g.f. has general term
$$u_n=\sum_{k=0}^{\frac{\lfloor n-1
\rfloor}{2}}\binom{n-1}{2k}C(k)\alpha^{n-2k-1}\beta^k.$$ In this
note, we shall be interested in the termwise ratios of the Hankel
transforms of the three sequences $u_n$, $u_n^*=u_{n+1}$ and $u_n^{**}=u_{n+2}$.
\newline\newline We will take the case $\alpha=-3$ and $\beta=-5$ to illustrate our results.
Thus let $a_n$ be the sequence with o.g.f. $\frac{x}{1-3x-5x^2}$. Then $a_n$ begins $0,1,3,14,57,241,\ldots$
with $$a_n=\sum_{k=0}^{\frac{\lfloor n-1
\rfloor}{2}}\binom{n-k-1}{k}3^{n-2k-1}5^k.$$ The series reversion of $\frac{x}{1-3x-5x^2}$ is $\frac{\sqrt{1+6x+29x^2}-3x-1}{10x}$ which generates
the sequence $u_n$ which begins $0,1,-3,4,18,-139,357,\ldots$ where
$$u_n=\sum_{k=0}^{\frac{\lfloor n-1
\rfloor}{2}}\binom{n-1}{2k}C(k)(-3)^{n-2k-1}(-5)^k.$$ We now form the shifted sequences
$$u_n^*=u_{n+1}=\sum_{k=0}^{\frac{\lfloor n
\rfloor}{2}}\binom{n}{2k}C(k)\alpha^{n-2k}\beta^k$$ and
$$u_n^{**}=u_{n+2}=\sum_{k=0}^{\frac{\lfloor n+1
\rfloor}{2}}\binom{n+1}{2k}C(k)\alpha^{n-2k+1}\beta^k.$$ We now let $h_n$, $h_n^*$ and $h_n^{**}$, respectively, be the Hankel transforms
of these sequences. Numerically, we find that the following:
\begin{center}
\begin{tabular}{|c|c|}\hline
Sequence & Hankel transform \\\hline
$u_n$    & $0,-1,-15,1750,890625,-2353515625,\ldots$\\\hline
$u_n*$   & $1,-5,-125,15625,9765625,-30517578125,\ldots$\\\hline
$u_n^{**}$ & $-3, -70, 7125, 3765625, -9843750000, -129058837890625,\ldots$\\
\hline
\end{tabular}\end{center}
These results suggest that $h_n^*=(-5)^{\binom{n+1}{2}}$, and
$$ \frac{(-1)^{n+1}h_{n+1}}{h_n^*}=a_{n+1}$$ along with
$$ \frac{(-1)^{n+1}h_n^{**}}{h_n^*}=a_{n+2}.$$
Thus in this case we obtain
$$h_{n+1}=(-1)^{n+1}(-5)^{\binom{n+1}{2}}a_{n+1}=(-1)^{n+1}(-5)^{\binom{n+1}{2}}\sum_{k=0}^{\frac{\lfloor n
\rfloor}{2}}\binom{n-k}{k}3^{n-2k}5^k$$ from which we infer that
$$h_n=(-1)^{n}(-5)^{\binom{n}{2}}\sum_{k=0}^{\frac{\lfloor n-1
\rfloor}{2}}\binom{n-k-1}{k}3^{n-2k-1}5^k.$$
Similarly, we find
$$h_n^{**}=(-1)^{n+1}(-5)^{\binom{n+1}{2}}a_{n+2}=(-1)^{n+1}(-5)^{\binom{n+1}{2}}\sum_{k=0}^{\frac{\lfloor n+1
\rfloor}{2}}\binom{n-k+1}{k}3^{n-2k+1}5^k.$$
Summarizing, we thus have
\begin{center}
\begin{tabular}{|c|c|}\hline
Sequence & Hankel transform \\\hline
$u_n$    & $h_n=(-1)^{n}(-5)^{\binom{n}{2}}\sum_{k=0}^{\frac{\lfloor n-1
\rfloor}{2}}\binom{n-k-1}{k}3^{n-2k-1}5^k$\\\hline
$u_n*$   & $h_n^*=(-5)^{\binom{n+1}{2}}$\\\hline
$u_n^{**}$ & $h_n^{**}=(-1)^{n+1}(-5)^{\binom{n+1}{2}}\sum_{k=0}^{\frac{\lfloor n+1
\rfloor}{2}}\binom{n-k+1}{k}3^{n-2k+1}5^k$\\
\hline
\end{tabular}\end{center}

We note that we have been able to recover the sequence $a_n$ in this example. Since the reversion of
the reversion of a series is the original series, we can now posit the
\begin{conjecture}
 Let $u_n$ be the sequence
$$u_n=\sum_{k=0}^{\frac{\lfloor n-1
\rfloor}{2}}\binom{n-1}{2k}C(k)\alpha^{n-2k-1}\beta^k$$ with integer parameters
$\alpha$ and $\beta$, and o.g.f.
$$u(x)=\frac{1-\alpha{x}-\sqrt{1-2\alpha{x}+(\alpha^2-4\beta)x^2}}{2\beta{x}}.$$
Let $h_n$ be the Hankel transform of $u_n$, $h_n^*$ the Hankel transform of $u_{n+1}$, and
$h_n^{**}$ be the Hankel transform of $u_{n+2}$. Further, let $a_n$ be the sequence with o.g.f. the
series reversion of $u(x)$, with
$$a_n=\sum_{k=0}^{\frac{\lfloor n-1
\rfloor}{2}}\binom{n-k-1}{k}\alpha^{n-2k-1}\beta^k.$$ Then
\begin{enumerate}
\item
$h_n^*=\beta^{\binom{n+1}{2}}$
\item
$ \frac{(-1)^{n+1}h_{n+1}}{h_n^*}=a_{n+1} \Rightarrow
h_n=(-1)^{n}\beta^{\binom{n}{2}}\sum_{k=0}^{\frac{\lfloor n-1
\rfloor}{2}}\binom{n-k-1}{k}\alpha^{n-2k-1}\beta^k$
\item
$ \frac{(-1)^{n+1}h_n^{**}}{h_n^*}=a_{n+2} \Rightarrow
h_n^{**}=(-1)^{n+1}\beta^{\binom{n+1}{2}}\sum_{k=0}^{\frac{\lfloor
n+1 \rfloor}{2}}\binom{n-k+1}{k}\alpha^{n-2k+1}\beta^k \qquad \qquad \Diamond$
\end{enumerate}
\end{conjecture}

Note that the Hankel transform of $u_{n+1}$, $h_n^*$ is independent of $\alpha$. This is due to
\begin{enumerate}
\item The binomial transform does not change the Hankel transform, and
\item The binomial transform of
$$\frac{1-\alpha{x}-\sqrt{1-2\alpha{x}+(\alpha^2-4\beta)x^2}}{2\beta x^2}$$ is given by
$$\frac{1-(\alpha+1){x}-\sqrt{1-2(\alpha+1){x}+((\alpha+1)^2-4\beta)x^2}}{2\beta x^2}.$$
\end{enumerate}
Elements in this family are related to coloured Motzkin paths. For other links between lattice paths and
Hankel transforms, see \cite{WW}.

The recovery of the sequence $a_n$ is an interesting feature of this family of sequences. That this is not always the case
is illustrated by the next example.
\end{example}
\begin{example} The family $\frac{x(1-\alpha x)}{1-\beta
x}$ for $\beta \ne 0$.
\newline\newline
The sequence with o.g.f. $\frac{x(1-\alpha x)}{1-\beta
x}$ has general term given by
$$a_n=(\beta-\alpha)\beta^{n-1}+\frac{\alpha}{\beta}0^n.$$ Here, $0^n$ is used to denote the  sequence
beginning $1,0,0,0,\ldots$ with o.g.f. $1$.
The reversion of the series $\frac{x(1-\alpha x)}{1-\beta
x}$, that is, the solution $u=u(x)$ of the equation
$$\frac{u(1-\alpha u)}{1-\beta
u}=x$$ is given by $$u(x)=\frac{1+\beta x-\sqrt{1-(2\alpha-\beta)2x+\beta^2x^2}}{2\alpha}.$$\end{example}
The sequence $u_n$ with this o.g.f. has general term
$$u_n=\sum_{k=0}^{n-1}\binom{n+k-1}{2k}C(k)\alpha^k(-\beta)^{n-k-1}.$$
Again, we shall be interested in the term-wise ratios of the Hankel transforms $h_n$, $h_n^*$ and $h_n^{**}$ respectively of the sequences
$u_n$, $u_n^*=u_{n+1}$ and $u_n^{**}=u_{n+2}$.
We obtain
\begin{conjecture} Using the notation above, we have
\begin{enumerate}
\item $h_n^*=(\alpha(\alpha-\beta))^{\binom{n+1}{2}}$
\item $h_{n+1}/h_n^*=\frac{(\alpha-\beta)^{n+1}-\alpha^{n+1}}{\beta} \Rightarrow h_n=\frac{(\alpha-\beta)^n-\alpha^n}{\beta}(\alpha(\alpha-\beta))^{\binom{n}{2}}$
\item $h_n^{**}/h_n^*=(\alpha-\beta)^{n+1} \Rightarrow h_n^{**}=(\alpha-\beta)^{n+1}(\alpha(\alpha-\beta))^{\binom{n+1}{2}} \qquad \qquad \Diamond$\end{enumerate}\end{conjecture}
We note that $h_{n+1}/h_n^*$ is the general term of the sequence with o.g.f. $\frac{-1}{(1-\alpha x)(1-(\alpha-\beta)x)}$ while
$h_n^{**}/h_n^*$ is the general term of the power sequence with o.g.f. $\frac{\alpha-\beta}{1-(\alpha-\beta)x}$. Thus in this case we do not recover terms of the
sequence with o.g.f. $\frac{x(1-\alpha x)}{1-\beta
x}$.
\begin{example} The family $x(1-\alpha x)$.
\newline\newline We note that this is in fact the case of $\frac{x(1-\alpha x)}{1-\beta x}$ where
$\beta=0$. The sequence with o.g.f. $x(1-\alpha x)$ is the sequence $0,1,-\alpha,0,0,0,\ldots$. The
reversion of the series $x(1-\alpha x)$, that is, the solution $u=u(x)$ of the equation
$$u(1-\alpha u)=x$$ is given by
$$u(x)=\frac{1-\sqrt{1-4\alpha x}}{2\alpha}.$$ For instance, the case $\alpha=1$ is that of the Catalan numbers
preceded by $0$. In general, $\frac{1-\sqrt{1-4\alpha x}}{2\alpha}$ is the o.g.f. of the sequence
$0,1,\alpha,2\alpha^2,5\alpha^3,14\alpha^4,\ldots$ with general term
$$ a_0=0, \qquad a_n=C(n-1)\alpha^{n-1}, n>0.$$
We obtain
\begin{conjecture} Using the notation above, we have
\begin{enumerate}
\item $h_n=-n\alpha^{n^2-1}, \qquad \frac{h_{n+1}}{h_n^*}=-(n+1)\alpha^n$
\item $h_n^*=\alpha^{n(n+1)}$
\item $h_n^{**}=\alpha^{(n+1)^2}, \qquad \frac{h_n^{**}}{h_n^*}=\alpha^{n+1}  \qquad \qquad \Diamond$
\end{enumerate}
\end{conjecture}
 We note in particular that this generalizes the well-known result
on the Hankel transforms of $C(n)$ and $C(n+1)$. We can in fact
easily modify the proof of the fact that the Hankel transform of
$C(n)$ is the all $1$'s sequence given in \cite{Rad1} to yield
\begin{proposition} The Hankel transform $h_n^*$ of the sequence
$C(n)\alpha^n$ is given by $h_n^*=\alpha^{n(n+1)}$.
\end{proposition}
\begin{proof} The coefficient of $x^{i+j+1}$ in $(1-\alpha
x)^2(1+\alpha x)^{2i+2j}$ is given by
$$\left\{\binom{2i+2j}{i+j+1}-2\binom{2i+2j}{i+j}+\binom{2i+2j}{i+j-1}\right\}\alpha^{i+j+1}=-2C(i+j)\alpha^{i+j+1}.$$
On the other hand, the coefficient of $x^k$ in $(1-\alpha
x)(1+\alpha x)^{2i}$ is equal to
$$\left\{\binom{2i}{k}-\binom{2i}{k-1}\right\}\alpha^k=\binom{2i}{k}\frac{2i-2k+1}{2i-k+1}\alpha^k.$$
Proceeding as in \cite{Rad1}, we obtain that
$$C(i+j)\alpha^{i+j}=\sum_{k=0}^{\min(i,j)}T_{i,k}T_{j,k}$$ where
$$T_{n,k}=\frac{\binom{2n}{n+k}(2k+1)}{n+k+1}\alpha^{n}.$$
Now $H_n=T_nU_n$ where $U_n$ is the transpose of $T_n$, where $T_n$
is the $(n+1) \times (n+1)$ matrix with general term $T_{i,k}$
\begin{displaymath}\mathbf{T_n}=\left(\begin{array}{ccccccc} 1 & 0 & 0
& 0 & 0 & 0 & \ldots \\\alpha & \alpha & 0 & 0 & 0 & 0 & \ldots \\
2\alpha^2 & 3\alpha^2 & \alpha^2 & 0 & 0 & 0 & \ldots \\ 5\alpha^3 & 9\alpha^3 & 5\alpha^3 & \alpha^3 & 0 & 0 & \ldots \\
14\alpha^4 & 28\alpha^4 & 20\alpha^4 & 7\alpha^4 & \alpha^4 & 0 &
\ldots
\\42\alpha^5 & 90\alpha^5 & 75\alpha^5 & 35\alpha^5 & 9\alpha^5 & \alpha^5 &\ldots\\ \vdots & \vdots & \vdots & \vdots & \vdots
& \vdots & \ddots\end{array}\right)\end{displaymath} Hence $h_n^*$
is the square of the product of the diagonal elements, that is
$$h_n^*=(\alpha^{\binom{n+1}{2}})^2=\alpha^{n(n+1)}.$$
\end{proof}
\end{example}

\bigskip
\hrule
\bigskip
\noindent 2000 {\it Mathematics Subject Classification}: Primary
11B83; Secondary 11C20, 11Y55

\noindent \emph{Keywords:} Hankel transform,
Catalan numbers, series reversion.


\begin{thebibliography}{9}

\bibitem{PasTri} P. Barry, On Integer-Sequence-Based Constructions of
Generalized Pascal Triangles, Journal of Integer Sequences, Vol 9, Article 06.2.4

\bibitem{Hankel1} A. Cvetkovic, P. Rajkovic and M. Ivkovic, Catalan Numbers, the Hankel Transform,
and Fiboncci Numbers, Journal of Integer Sequences, 5, May 2002,
Art. 02.1.3

\bibitem{Krat} C. Krattenthaler, Advanced Determinant Calculus, available electronically at
\texttt{http://www.mat.univie.ac.at/~slc/wpapers/s42kratt.pdf},2007

\bibitem{Layman} J. W. Layman, The Hankel Transform and some of its
properties, Journal of Integer Sequences, Vol. 4, (2001), Article
01.1.5

\bibitem{PPWW} P. Peart, W-J. Woan, Generating Functions via Hankel and Stieltjes
Matrices, Journal of Integer Sequences, Vol. 3 (2000), Article
00.2.1

\bibitem{Hankel2}P. Rajkovic, M. D. Petkovic and P. Barry, The Hankel Transform of the Sum of Consecutive Generalized Catalan Numbers,
available electronically at
\texttt{http://www.arXiv.org/math/pdf/0604/0604422.pdf},2007.

\bibitem{Rad1} C. Radoux, D\'eterminant de Hankel construit sur des p\^olynomes li\'es aux nombres de d\'erangements, European J. Combin., 12 (1991) 327–-329.

\bibitem{Rad2} C. Radoux, Nombres de Catalan g\'en\'eralis\'es, Bull.
Belg. Math. Soc. 4(1997), 289--292

\bibitem{Rad3} C. Radoux, Calcul effectif de certains d\'eterminants de Hankel, Bull. Soc. Belg. S\'er. B 31 (1) (1979) 49--55.

\bibitem{Rad4} C. Radoux, D\'eterminant de Hankel construit sur les p\^olynomes de Hermite, Ann. Soc. Sci. Bruxelles 104 (2) (1990) 59--61

\bibitem{Rad5} C. Radoux, Addition formulas for polynomials built on classical combinatorial sequences, J. of Comp. and App. Math. 115 (2000) 471--477

\bibitem{WW} W-J. Woan, Hankel Matrices and Lattice Paths, Journal
of Integer Sequences, Vol. 4 (2001), Article 01.1.2

\bibitem{Zeil} D. Zeilberger, An Umbral Approach to the Hankel Transform for Sequences, available electronically at
\texttt{http://www.math.rutgers.edu/~zeilberg/mamarim/mamarimPDF/hankel.pdf},2007.

\end{thebibliography}
\end{document}